\documentclass{aims}
\usepackage{amsmath}
  \usepackage{paralist}
  \usepackage{graphics} 
  \usepackage{epsfig} 
\usepackage{graphicx}  \usepackage{epstopdf}
 \usepackage[colorlinks=true]{hyperref}
 \hypersetup{urlcolor=blue, citecolor=red}
\usepackage{wrapfig}
 \usepackage{amsfonts}
\usepackage{amssymb}
\hypersetup{urlcolor=blue, citecolor=red}
\usepackage{hyperref}
\usepackage[all]{xy}

\def\currentvolume{X}
 \def\currentissue{X}
  \def\currentyear{200X}
   \def\currentmonth{XX}
    \def\ppages{X--XX}

\theoremstyle{definition}

\newcommand{\eps}{\varepsilon}
\numberwithin{equation}{section}

\title[Classifying Homogenization-Based Numerical Mathods]
      {An Attempt at Classifying Homogenization-Based Numerical Mathods}

\author[Emmanuel Fr\'enod]{}


 \email{emmanuel.frenod@univ-ubs.fr}

\newcommand{\opEquLim}{{O}}
\newcommand{\opEquOsc}{{O^{\eps}}}
\newcommand{\opEquLimDis}{{O_{\Delta z}}}
\newcommand{\opEquLimDisUN}{{O_{\Delta z}^1}}
\newcommand{\opEquLimDisOsc}{{O_{\Delta z}^{\eps}}}
\newcommand{\opTscLim}{{\mathcal{O}}} 
\newcommand{\opTscOsc}{{\mathcal{O}^{\eps}}}
\newcommand{\opTscLimDis}{{\mathcal{O}_{\Delta z}}}
\newcommand{\opTscLimDisUN}{{\mathcal{O}_{\Delta z}^1}}
\newcommand{\opTscLimDisOsc}{{\mathcal{O}_{\Delta z}^{\eps}}}

\newcommand{\OpMacro}{{M}}

\newcommand{\OpMacroEpsDis}{{M^\eps_{\Delta z}}}
\newcommand{\OpMacroDis}{{M_{\Delta z}}}

\begin{document}
\maketitle

\centerline{\scshape Emmanuel Fr\'enod}
\medskip
{\footnotesize
 \centerline{Universit\'e de Bretagne-Sud,  UMR 6205, LMBA, F-56000 Vannes, France}
   \center  \centerline{AND}
   \centerline{Projet INRIA Calvi, Universit\'{e} de Strasbourg, IRMA,}
   \centerline{7 rue Ren\'e Descartes, F-67084 Strasbourg Cedex, France}
} 

\medskip
 \centerline{(Communicated by the associate editor name)}

\pagestyle{myheadings}
 \renewcommand{\sectionmark}[1]{\markboth{#1}{}}
\renewcommand{\sectionmark}[1]{\markright{\thesection\ #1}}
%

\begin{abstract} 
In this note, a classification of Homogenization-Based Numerical Methods and (in particular) of Numerical Methods
that are based on the Two-Scale Convergence is done. In this classification stand: Direct Homogenization-Based Numerical Methods;
H-Measure-Based Numerical Methods; Two-Scale Numerical Methods and TSAPS: Two-Scale Asymptotic Preserving Schemes.
\end{abstract}

\section{Introduction}
A Homogenization-Based Numerical Method is a numerical method that incorporates in its conception concepts coming from Homogenization
Theory. Doing this gives to the built method the capability to tackle efficiently heterogeneities or oscillations. 
This approach can be applied to problems occurring in a heterogeneous medium, that have oscillating boundary conditions or that are
constrained to oscillate by an external action (for instance a magnetic field on a charged particle cloud).

This topic is currently active. The goal of this special issue is to emphasis recent advances in this topic in a wide variety of application
fields.
\\

This introductory paper introduces a classification of Homogenization-Based Numerical Methods, in which  stand: 
Direct Homogenization-Based Numerical Methods;
H-Measure-Based Numerical Methods; Two-Scale Numerical Methods and TSAPS: Two-Scale Asymptotic Preserving Schemes.

\section{Direct Homogenization-Based Numerical Methods}
The context of  Direct Homogeniza\-tion-Based Numerical Methods is depicted in the next diagram:
\label{DHBNM} 
\begin{equation} 
\label{Diag-DHBNM} 
 \xymatrix{
 \textnormal{
   \begin{tabular}{c}
   {\small$u^{\eps}$ solution to \vspace{3pt}}\\
   {\large$\opEquOsc\,u^{\eps} = 0$}
   \end{tabular}
 }
  \ar[rr]_{\eps \, \to \, 0}
 & &
 \textnormal{
   \begin{tabular}{c}
   {\small$u$ solution to \vspace{3pt}}\\
   {\large$\opEquLim\,u = 0$}
    \end{tabular}
  }  
  \\
  & & \\
 ~
  & &
 \textnormal{
  \begin{tabular}{c}
  {\small$u_{\Delta z}$ solution to \vspace{3pt}}\\
  {\large$\opEquLimDis\,u_{\Delta z} = 0$}
  \end{tabular}
  \ar[uu]^{\Delta z \, \to \, 0}
 }
 }
\end{equation}
It is when we face with an operator $\opEquOsc$ that generates
in solution $u^{\eps}$ of equation $\opEquOsc\,u^{\eps} = 0$ oscillations or heterogeneities of characteristic 
size $\eps$ -\,which is small\,- and when it is 
known that, in some sense, for small $\eps$, $u^{\eps}(z)$ is close to $u(z)$  for which is known a well-posed problem $\opEquLim\,u = 0$.

In this context, it is possible, in place of building a numerical approximation of operator $\opEquOsc$, to build a numerical
operator $\opEquLimDis$ approximating $\opEquLim$. Then solving  $\opEquLimDis\,u_{\Delta z}$ gives a solution $u_{\Delta z}(z)$ 
which is close to $u$
and consequently to $u^{\eps}(z)$, when $\eps$ is small. 
This approach permits to obtain an approximation of $u^{\eps}(z)$ without resolving the oscillations 
the model to compute it contains.
\\

In the case when a corrector result is known, {\it i.e}, if in association with $u(z)$, a function $u^1(z)$, solution to well-posed 
equation $\opEquLim^1\,u^1 = 0$, is such that $u^{\eps}(z)$ is close to $u(z)+\eps u^1(z)$ for small $\eps$, it is possible 
build two numerical operators $\opEquLimDis$ and $\opEquLimDisUN$ that are discretizations of $\opEquLim$ and $\opEquLim^1$.
Using them, we can compute approximated solutions $u_{\Delta z}(z)$ and $u_{\Delta z}^1(z)$ of $u(z)$ and $u^1(z)$ and obtain
a good approximation of $u^{\eps}(z)$ computing $u_{\Delta z}(z)+\eps u_{\Delta z}^1(z)$. Such a method is called order-1 
Direct Homogenization-Based Numerical Methods and is illutrated by the following diagram.
\begin{equation} 
\label{Diag-DHBNM-2} 
 \xymatrix{
 \textnormal{
   \begin{tabular}{c}
   {\small$u^{\eps}$ solution to \vspace{3pt}}\\
   {\large$\opEquOsc\,u^{\eps} = 0$}
   \end{tabular}
 }
  \ar[rr]_{\eps \, \to \, 0}
 & &
 \textnormal{
   \begin{tabular}{c}
   {\small$u$, $u_1$ solutions to \vspace{3pt}}\\
   {\large$\opEquLim\,u = 0$\vspace{3pt}}\\
   {\large$\opEquLim^1\,u^1 = 0$}
    \end{tabular}
  }  
  \\
  & & \\
 ~
  & &
 \textnormal{
  \begin{tabular}{c}
  {\small$u_{\Delta z}$, \small$u_{\Delta z}$ solution to \vspace{3pt}}\\
  {\large$\opEquLimDis\,u_{\Delta z} = 0$\vspace{3pt}}\\
  {\large$\opEquLimDisUN\,u_{\Delta z}^1 = 0$}
  \end{tabular}
  \ar[uu]^{\Delta z \, \to \, 0}
 }
 }
\end{equation}
\\

The paper by Legoll \& Minvielle  \cite{LegMin}, by Laptev \cite{Laptev}, by Bernard, Fr\'enod \& Rousseau  \cite{BerFreRou2},
and by Xu \& Yue \cite{XuYue} of this special issue may enter this framework
%
\section{H-Measure-Based Numerical Methods} 
The context of those kind of methods is when the transition from a microscopic scale (of size $\eps$) 
to a macroscopic one (of size $1$) -\,with quantities of interest at the microscopic scale that
are not the same as the quantities of interest at the macroscopic scale\,- needs to be described. 
This occurs for instance in the simulation of phenomena some parts of which call upon quantum description
or in the simulation of turbulence. 
This context can be represented by the following diagram:
\label{HMBNM} 
\begin{equation} 
\label{Diag-AHMBNM} 
 \xymatrix{
 \textnormal{
   \begin{tabular}{c}
   {\small$u^{\eps}$ solution to \vspace{3pt}}\\
   {\large$\opEquOsc\,u^{\eps} = 0$} \\
   {$e^\eps = \mathbf{E}(u^{\eps})$}
   \end{tabular}
 }
  \ar[rr]_{\eps \, \to \, 0}
 & &
 \textnormal{
   \begin{tabular}{c}
   {\small$e$ solution to \vspace{3pt}}\\
   {\large$\OpMacro\,e = 0$}
    \end{tabular}
  }  
  \\
  & & \\
 \textnormal{
  \begin{tabular}{c}
  {\small$u_{\Delta z}^{\eps}$ solution to \vspace{3pt}}\\
  {\large$\OpMacroEpsDis\,e_{\Delta z}^{\eps} = 0$}
  \end{tabular}
 }
 \ar[uu]^{\Delta z \, \to \, 0}  \ar[rr]_{\eps \, \to \, 0}
  & &
 \textnormal{
  \begin{tabular}{c}
  {\small$u_{\Delta z}$ solution to \vspace{3pt}}\\
  {\large$\OpMacroDis\,e_{\Delta z} = 0$}
  \end{tabular}
  \ar[uu]^{\Delta z \, \to \, 0}
 }
 }
\end{equation}
and explained as follows.
The part in the top left of diagram \eqref{Diag-AHMBNM} symbolizes a problem which is set at the microscopic level.
This problem writes $ \opEquOsc\,u^{\eps} = 0$ and generates oscillations in its solution $u^{\eps}$.
Besides, the quantity that makes sense at the macroscopic level is $e^\eps$; it is related to $u^{\eps}$ by
a non-linear relation $e^\eps = \mathbf{E}(u^{\eps})$ and it is, in some sense, close to $e$ solution to $\OpMacro\,e = 0$
(see the top right of the diagram) which represents the model at the macroscopic level.

Then,  the goal of a H-Measure-Based Numerical Method consists in building a numerical operator $\OpMacroEpsDis$, 
giving a numerical solution $e_{\Delta z}^{\eps} $ close to $e^\eps$, for any $\eps$ as soon as $\Delta z$ is small
(see the bottom left of the diagram), which 
behaves as a numerical approximation of $M$ when $\eps$ is small (see the bottom right of the diagram).
\\

The paper by Tartar \cite{Tartar2013} of this special issue lays the foundation of the theory for those kinds of methods.
\section{Two-Scale Numerical Methods}
The papers by Assyr, Bai \& Vilmar  \cite{AbdBaiVil},
Back \& Fr\'enod  \cite{BackFre},
Faye, Fr\'enod \& Seck \cite{FayFreSec2013},
Fr\'enod, Hirtoaga \& Sonnendr\"ucker \cite{FreHirSon},
Lutz \cite{Lutz2013} and 
Henning \& Ohlberger  \cite{HenOhl} of this special issue  are related to this framework of Two-Scale Numerical Methods.
\\

An order-0 Two-Scale Numerical Method may be explained using the following diagram:
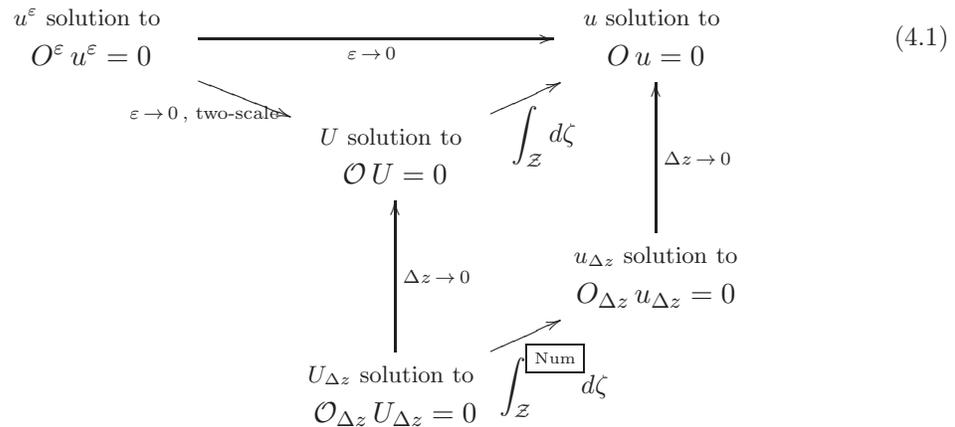
\label{TSNM} 
\begin{equation}
\label{Diag-TSNM-1} 
\xymatrix@C=1pc@R=1pc{
\textnormal{
\begin{tabular}{c}
{\small $u^{\eps}$ solution to  \vspace{3pt}} \\ {\large$\opEquOsc\,u^{\eps} = 0$} 
\end{tabular}
}
\ar[rrr]_{\eps \, \to \, 0 \,} \ar[drr]_{\eps \, \to\,0\,, \textnormal{ two-scale}~~} 
& & & 
\textnormal{
\begin{tabular}{c}
{\small $u$ solution to \vspace{3pt}} \\ {\large $\opEquLim\,u=0$}
\end{tabular}
} \\
& & \textnormal{
\begin{tabular}{c}
{\small $U$ solution to \vspace{3pt}} \\ {\large$\opTscLim\,U = 0$}
\end{tabular}
} \ar[ru]_{\displaystyle \int_{\mathcal{Z}}d\zeta} & \\
\textnormal{
} & & & \textnormal{
\begin{tabular}{c}
{\small $u_{\Delta z}$ solution to\vspace{3pt}} \\ {\large $\opEquLimDis\,u_{\Delta z} = 0$}
\end{tabular}
} \ar[uu]_{\Delta z \, \to\,0} \\
& & \textnormal{
\begin{tabular}{c}
{\small $U_{\Delta z}$ solution to  \vspace{3pt}} \\{\large $\opTscLimDis\,U_{\Delta z} = 0$}
\end{tabular}
} \ar[uu]_{\Delta z \, \to\, 0} \ar[ru]_{\displaystyle \int_{\mathcal{Z}}^{\fbox{\tiny Num}}\!\!d\zeta} & \\
}
\end{equation} 
The context includes the one of Direct Homogenization-Based Numerical Methods and diagram  \eqref{Diag-TSNM-1} has to be 
regarded as a prism. Its deepest layer is nothing but diagram \eqref{Diag-DHBNM}. Yet, if more is known about the asymptotic behavior
of $u^{\eps}$, {\it i.e.} if it is known that $u^\eps(z)$ is close to $U(z,\frac{z}{\eps})$, with $U(z,\zeta)$ periodic in $\zeta$, when
$\eps$ is small (which can be translated as  $u^\eps(z)$ Two-Scale Converges to $U(z,\zeta)$) and if it is known a well posed
problem $\opTscLim\,U = 0$ for $U$ (see the middle of the diagram), that gives the equation for $u$ (see the top right of the diagram)
when integrated with respect to periodic variable $\zeta$, it is possible to build a specific numerical method.

This method consists in building a numerical approximation $\opTscLimDis$ of operator $\opTscLim$. Using this operator 
can give a numerical solution $U_{\Delta z}$ (see the  bottom of the diagram) and $U_{\Delta z}(z,\frac{z}{\eps})$ is an approximation
of $u^\eps(z)$ for small $\eps$. To be consistent with the continuous level, a numerical integration of the $\opTscLimDis\,U_{\Delta z} = 0$ needs to
yield a numerical approximation of the equation for $u$ (see the bottom right of the diagram).
\\

When a little more is known concerning the asymptotic behavior of $u^\eps$ when $\eps$ is small, {\it i.e.} if $u^\eps$ is close to
$U(z,\frac{z}{\eps})+\eps U^1(z,\frac{z}{\eps})$
with  $U^1(z,\zeta)$ also periodic in $\zeta$ and if a well-posed problem is known for $U^1$, we can enrich diagram \eqref{Diag-TSNM-1}
and obtain the following diagram of order-1 Two-Scale Numerical Methods:
\begin{equation}
\label{Diag-TSNM-2} 
\xymatrix@C=1pc@R=1pc{
\textnormal{
\begin{tabular}{c}
{\small $u^{\eps}$ solution to  \vspace{3pt}} \\ {\large$\opEquOsc\,u^{\eps} = 0$}
\end{tabular}
}
\ar[rrr]_{\eps \, \to \, 0 } \ar[drr]_{\eps \, \to\,0\,, \textnormal{ two-scale}~~} 
& & & 
\textnormal{
\begin{tabular}{c}
{\small $u$, $u_1$ solutions to \vspace{3pt}} \\ {\large $\opEquLim\,u=0$} \\ {\large$\opEquLim^1u^1=0$}
\end{tabular}
} \\
& & \textnormal{
\begin{tabular}{c}
{\small $U$, $U^1$ solutions to \vspace{3pt}} \\ {\large$\opTscLim\,U = 0$} \\ {\large$\opTscLim^1U^1 = 0$}
\end{tabular}
} \ar[ru]_{\displaystyle \int_{\mathcal{Z}}d\zeta} & \\
\textnormal{
} & & & \textnormal{
\begin{tabular}{c}
{\small $u_{\Delta z}$, $u^1_{\Delta z}$ solutions to\vspace{3pt}} \\ 
 {\large $\opEquLimDis\,u_{\Delta z} = 0$} \\ {\large$\opEquLimDisUN u^1_{\Delta z} = 0$}
\end{tabular}
} \ar[uu]_{\Delta z \, \to\,0} \\
& & \textnormal{
\begin{tabular}{c}
{\small $U_{\Delta z}$, $U^1_{\Delta z}$  solutions to  \vspace{3pt}} \\
  {\large $\opTscLimDis\,U_{\Delta z} = 0$} \\ {\large $\opTscLimDisUN U^1_{\Delta z} = 0$}
\end{tabular}
} \ar[uu]_{\Delta z \, \to\, 0} \ar[ru]_{\displaystyle \int_{\mathcal{Z}}^{\fbox{\tiny Num}}\!\!d\zeta} & \\
}
\end{equation} 
%
\section{TSAPS: Two-Scale Asymptotic Preserving Schemes}
\label{TSAPS} 
To describe Two-Scale Asymptotic Preserving Schemes, it is first needed to describe what is an
Asymptotic Preserving Scheme (or AP-Scheme in short).
\begin{equation} 
\label{Diag-AP} 
 \xymatrix{
 \textnormal{
   \begin{tabular}{c}
   {\small$u^{\eps}$ solution to \vspace{3pt}}\\
   {\large$\opEquOsc\,u^{\eps} = 0$}
   \end{tabular}
 }
  \ar[rr]_{\eps \, \to \, 0}
 & &
 \textnormal{
   \begin{tabular}{c}
   {\small$u$ solution to \vspace{3pt}}\\
   {\large$\opEquLim\,u = 0$}
    \end{tabular}
  }  
  \\
  & & \\
 \textnormal{
  \begin{tabular}{c}
  {\small$u_{\Delta z}^{\eps}$ solution to \vspace{3pt}}\\
  {\large$\opEquLimDisOsc\,u_{\Delta z}^{\eps} = 0$}
  \end{tabular}
 }
 \ar[uu]^{\Delta z \, \to \, 0}  \ar[rr]_{\eps \, \to \, 0}
  & &
 \textnormal{
  \begin{tabular}{c}
  {\small$u_{\Delta z}$ solution to \vspace{3pt}}\\
  {\large$\opEquLimDis\,u_{\Delta z} = 0$}
  \end{tabular}
  \ar[uu]^{\Delta z \, \to \, 0}
 }
 }
\end{equation}
For this, we comment on diagram  \eqref{Diag-AP}.
The context is when we are face-to-face with an operator $\opEquOsc$ which is approached, when $\eps$ is small, 
by another operator $\opEquLim$ which has not the same nature as $\opEquOsc$.
An Asymptotic Preserving Scheme to approximate problem $\opEquOsc\,u^{\eps} = 0$ (see the top left of the diagram) 
is a numerical operator $\opEquLimDisOsc$ that gives, when solving $\opEquLimDisOsc\,u_{\Delta z}^{\eps} = 0$ 
(see the bottom right of the diagram) a numerical solution
$u_{\Delta z}^{\eps}$ which is close to $u$, with an accuracy depending on step $\Delta z$ and not on $\eps$.
Besides, this operator needs to mimic the behavior of an numerical approximation (see the bottom right of the diagram) 
of limit problem $\opEquLim\,u = 0$ (see the top right of the diagram) as $\eps$ is small.

For an introduction to this kind of method the reader is referred to {Jin \cite{Jin:SIAM:1999}}.
\\

The explanation of TSAPS, will be based on the following diagram:
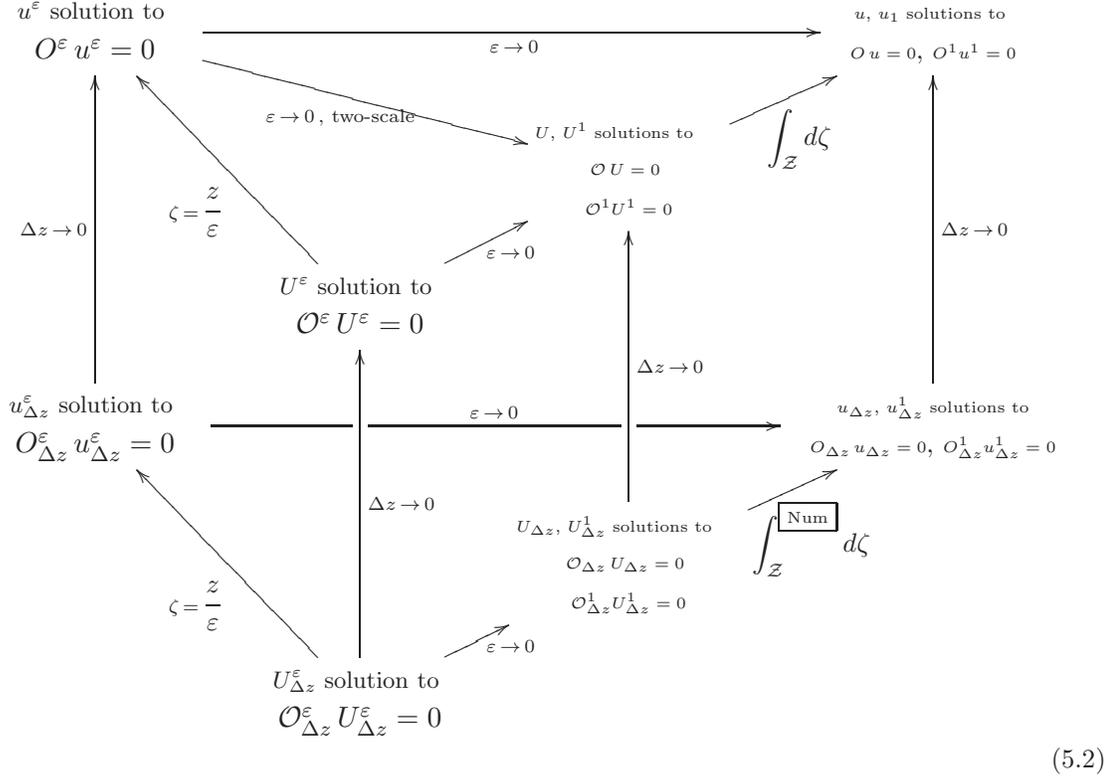
{
\begin{equation}
\label{Diag-TSAPS} 
\xymatrix@C=1pc@R=1pc{
\textnormal{
\begin{tabular}{c}
{\small $u^{\eps}$ solution to  \vspace{3pt}} \\ {\large$\opEquOsc\,u^{\eps} = 0$}
\end{tabular}
}
\ar[rrr]_{\eps \, \to \, 0 } \ar[drr]_{\eps \, \to\,0\,, \textnormal{ two-scale}} 
& & & 
\textnormal{
\begin{tabular}{c}
{\tiny $u$, $u_1$ solutions to \vspace{3pt}} \\ {\tiny $\opEquLim\,u=0$}, {\tiny$\opEquLim^1u^1=0$}
\end{tabular}
} \\
& & \textnormal{
\begin{tabular}{c}
~\hspace{-15pt}
{\tiny $U$, $U^1$ solutions to \vspace{2pt}} \\ {\tiny$\opTscLim\,U = 0$ \vspace{3pt}} \\  {\tiny$\opTscLim^1U^1 = 0$}
\end{tabular}
} \ar[ru]_{\displaystyle \int_{\mathcal{Z}}d\zeta} & \\
& \textnormal{
\begin{tabular}{c}
{\small $U^{\eps}$ solution to  \vspace{3pt}}\\ {\large$\opTscOsc\,U^{\eps} = 0$}
\end{tabular}
} \ar[luu]^{\zeta\,=\,\cfrac{z}{\eps}} \ar[ru]_{\eps\,\to\,0} & & \\
\textnormal{
\begin{tabular}{c}
{\small $u_{\Delta z}^{\eps}$ solution to \vspace{3pt} }\\ {\large$\opEquLimDisOsc \,u_{\Delta z}^{\eps} = 0$}
\end{tabular}
} \ar[uuu]^{\Delta z \, \to\, 0} \ar@{->}'[r]
'[rr]^{\eps\,\to\,0}
[rrr] & & & \textnormal{
\begin{tabular}{c}
{\tiny $u_{\Delta z}$, $u^1_{\Delta z}$ solutions to\vspace{3pt}} \\ 
 {\tiny $\opEquLimDis\,u_{\Delta z} = 0$},  {\tiny $\opEquLimDisUN u^1_{\Delta z} = 0$}
\end{tabular}
} \ar[uuu]_{\Delta z \, \to\,0} \\
& & \textnormal{
\begin{tabular}{c}
~\hspace{-15pt}
{\tiny $U_{\Delta z}$, $U^1_{\Delta z}$  solutions to  \vspace{2pt}} \\
  {\tiny $\opTscLimDis\,U_{\Delta z} = 0$ \vspace{3pt}} \\ {\tiny $\opTscLimDisUN U^1_{\Delta z} = 0$}
\end{tabular}
} \ar[uuu]_{\Delta z \, \to\, 0} \ar[ru]_{\displaystyle \int_{\mathcal{Z}}^{\fbox{\tiny Num}}d\zeta} & \\
& \textnormal{
\begin{tabular}{c}
{\small $U_{\Delta z}^{\eps}$ solution to \vspace{3pt} }\\ {\large$\opTscLimDisOsc \, U_{\Delta z}^{\eps} = 0$}
\end{tabular}
} \ar[luu]^{\zeta\,=\,\cfrac{z}{\eps}} \ar[uuu]_{\Delta z\,\to\,0} \ar[ru]_{\eps\,\to\,0} & & 
}
\end{equation} 
}
This diagram has to be regarded as a prism with three layers. The deepest one is the diagram of the AP-schemes.
At the top left of this layer is found the equation that generates in its solution oscillations of size $\eps$. 
At the top right stands the limit problem, as $\eps$ is small. 
(This limit problem is assumed to be an order-1 problem, {\it i.e.} $u^\eps \sim u^0 + \eps u^1$ for $\eps$ small and 
equations are known for $u^0$ and $u^1$.)
At the bottom left stands the AP-Scheme that approximate
equation $\opEquOsc\,u^{\eps} = 0$ for any $\eps$ and that mimics an approximation of the limit problem when
$\eps$ is small (see the bottom right of the layer).
\\
The middle layer is exactly the diagram of the order-1 Two-Scale Numerical Methods.
\\
The top layer is the new part. at the bottom, stands the TSAPS. This is a numerical method that gives a solution
$U_{\Delta z}^{\eps}$ which depends on two variables $z$ and $\zeta$. When taken in $\zeta=z/\eps$, $U_{\Delta z}^{\eps}$
gives a numerical approximation of the solution to the problem given at the top right of the diagram, with an accuracy that
only depends on the discretization step $\Delta z$, and not on $\eps$. Moreover, as $\eps$ is small, the TSAPS $\opTscLimDisOsc$
needs to mimic the behavior of the order-1 Two-Scale Numerical Operator (the couple ($\opTscLim$, $\opTscLim^1$)).
To builtd a TSAPS, a reformulation of problem $\opEquOsc\,u^{\eps} = 0$ calling upon a Two-Scale Macro-Micro Decomposition
(that reads $\opTscOsc\,U^{\eps} = 0$, see the middle of the diagram) is used.
A first step towards TSAPS is led in 
{{Crouseilles}, {Fr\'enod}, {Hirstoaga} \& {Mouton} \cite{crouseilles:hal-00638617}}.

\bibliographystyle{plain}
\bibliography{biblio}

\end{document}